\documentclass[12pt]{amsart}   
\pdfoutput=1 

\usepackage{graphicx,subfigure,color}
\usepackage{latexsym,amsmath,amsfonts,amscd,amsthm,epstopdf,lineno}
\usepackage[T1]{fontenc}
\usepackage[colorlinks=true,citecolor=blue,urlcolor=blue]{hyperref}

\usepackage{caption,here}
\usepackage{textcomp}
\usepackage{tikz}
\usetikzlibrary{decorations.pathmorphing,decorations.pathreplacing}
\usepackage{siunitx}
\usepackage{algorithm}  
\usepackage{enumitem}   

\usepackage{fullpage}

\newtheorem{theorem}{Theorem}[section]
\newtheorem{remark}{Remark}
\newtheorem{lemma}[theorem]{Lemma}

\title{Higher order A-stable schemes for the wave equation using a recursive convolution approach}

\author[M. Causley]{M. Causley}
\address{Department of Mathematics, Michigan State University, East Lansing, MI 48824}
\email{causleym@math.msu.edu}
\thanks{This work has been supported in part by AFOSR grants FA9550-11-1-0281, FA9550-12-1-0343 and FA9550-12-1-0455, NSF grant DMS-1115709, and MSU Foundation grant SPG-RG100059.}

\author[A. Christlieb]{A. Christlieb}

\begin{document}

\maketitle


\date{}

\begin{abstract}

In several recent works \cite{Causley2013a}, \cite{Causley2013}, we developed a new second order, A-stable approach to wave propagation problems based on the method of lines transpose (MOL$^T$) formulation combined with alternating direction implicit (ADI) schemes. Because our method is based on an integral solution of the ADI splitting of the MOL$^T$ formulation, we are able to easily embed non-Cartesian boundaries and include point sources  with exact spatial resolution. Further, we developed an efficient $O(N)$ convolution algorithm for rapid evaluation of the solution, which makes our method competitive with explicit finite difference (e.g., FDTD) solvers, both in terms of accuracy and time to solution, even for Courant numbers slightly larger than 1. We have demonstrated the utility of this method by applying it to a range of problems with complex geometry, including cavities with cusps.


In this work, we present several important modifications to our recently developed wave solver. We obtain a family of wave solvers which are unconditionally stable, accurate of order $2P$, and require $O(P^d N)$ operations per time step, where $N$ is the number of spatial points, and $d$ the number of spatial dimensions. We obtain these schemes by including higher derivatives of the solution, rather than increasing the number of time levels. The novel aspect of our approach is that the higher derivatives are constructed using successive applications of the convolution operator.

We develop these schemes in one spatial dimension, and then extend the results to higher dimensions, by reformulating the ADI scheme to include recursive convolution. Thus, we retain a fast, unconditionally stable scheme, which does not suffer from the large dispersion errors characteristic to the ADI method. We demonstrate the utility of the method by applying it to a host of wave propagation problems. This method holds great promise for developing higher order, parallelizable algorithms for solving hyperbolic PDEs, and can also be extended to parabolic PDEs.

\end{abstract}



\section{Introduction}

In recent works \cite{Causley2013a}, \cite{Causley2013}, the method of lines transpose (MOL$^T$) has been utilized to solve the wave equation, resulting in a second order accurate, A-stable numerical scheme. The solution is constructed using a boundary-corrected integral equation, which is derived in a semi-discrete setting. Thus, the solution at time $t^{n+1}$ is found by convolving the solution at previous time levels against the semi-discrete Green's function. Upon spatial discretization, traditional convolution requires $O(N^2)$ operations for a total of $N$ spatial points per time step. However, we have developed a fast convolution algorithm for the one-dimensional problem, which reduces the computational cost to  $O(N)$ operations \cite{Causley2013a}. This efficiency was additionally extended to the wave equation in higher dimensions by applying alternate direction implicit (ADI) splitting to the semi-discrete elliptic differential operator.

Our solver is intended to act as a field solver for Maxwell's equations in plasma simulations. In this regard, the MOL$^T$ approach has three distinct advantages:  it can capture time-dependent point sources (particles) with exact spatial resolution; it is $O(N)$, A-stable, and second order accurate; and it can incorporate complex geometries by embedding boundaries in a Cartesian mesh. However, we also point out that our methods are quite suitable for general electromagnetics, and acoustic problems, and are competitive to standard explicit finite difference methods (for example, the traditional FDTD scheme of Yee), both in terms of computational complexity and accuracy, but without the CFL time step restriction. 

While there is no stability restriction placed on our A-stable scheme, considerations of accuracy present themselves when the CFL number becomes large. Indeed, when large time steps are taken, the anisotropies introduced by the dimensional splitting are very pronounced. This problem has also been observed in the FDTD-ADI implementation of Maxwell's equations, which were introduced simultaneously by Namiki \cite{Namiki1999} and Zheng et. al. \cite{Zheng1999}. The splitting error can be understood as a dispersive term in the leading order of the truncation error \cite{Namiki2000}. Fornberg et. al. \cite{Fornberg2001,Fornberg2007} have studied the problem in great detail, and shown that if the dispersion error can be depressed with higher order spatial resolution, the resulting scheme for Maxwell's equations can be lifted to higher order in time using Richardson extrapolation, thus removing the ADI anisotropy.

In this work, rather than working with the first order Maxwell formulation, we shall apply ADI splitting directly to the wave equation. This is not a new idea, and in fact was first proposed by Lees \cite{Lees1962} shortly after Peaceman and Rachford applied it to the heat equation. Lees built upon the pioneering work of Von-Neumann, who first proposed an implicit finite difference solution for the 1d wave equation, which when viewed in a semi-discrete sense is essentially identical to our equation \eqref{eqn:Semi} below.

The notion of obtaining higher order ADI algorithms is also not novel, having first been presented by Fairweather \cite{Fairweather1965}. But the Fairweather's approach did not remain A-stable, as was the case for the second order scheme of Lees. The work on higher order ADI implementations for second order hyperbolic PDEs continues to this day, with recent emphasis placed on the use of compact finite differences \cite{Deng2012}. But what of higher order schemes which are also A-stable?

In the pioneering work by Dahlquist \cite{Dahlquist1963}, it is stated that no linear multistep scheme applied to the problem $y' = f(x,y)$ can simultaneously achieve A-stability, and order of accuracy greater than 2. Slightly less well know is that a decade later, the same result was proven again by Dahlquist \cite{Dahlquist1978} for periodic initial value problems of the form $y" = f(x,y)$. But the Dahlquist barrier can be broken, by not limiting the ODE solver to a linear multistep scheme, a fact pointed out by Ehle \cite{Ehle1968}, in his study of multistage implicit Runge-Kutta, as well as multiderivative schemes.

As such, a linear multiderivative scheme can achieve higher orders of accuracy, and remain A-stable. This result has been known for decades in the solution of periodic initial value problems, beginning with the work of Numerov, and Lambert and Watson \cite{Lambert1976}, and remains active to this day \cite{Stavroyiannis2009}, \cite{Stavroyiannis2010}.

On the other hand, multiderivative methods for solving hyperbolic PDEs has been considered much more sparsely in the literature. This could be attributed to the fact that, in contrast to ODEs, the introduction of spatial dependence creates several complications. In particular, we now must consider boundary conditions, and how the inclusion of higher derivatives effects the solution near the boundaries. Additionally, there is the issue of computational complexity. For instance, the traditional tridiagonal solves used in finite difference algorithms will now be replaced by banded matrices of growing bandwidth, which must be inverted at each time step. These issues were addressed very recently in \cite{Luo2013} for the telegraph equation, where an implicit Hermite interpolation in time was used to achieve a fourth order, A-stable numerical scheme. The 3-point spatial stencil was maintained, by using fourth order compact finite differences, and consistent endpoint corrections were derived for Dirichlet boundary conditions.

In addition to enhancing the stability properties of higher order methods, the use of multiderivative schemes also holds great promise for computational efficiency in parallel codes. This case was made recently in \cite{Seal2013}, where multiderivative methods are developed for hyperbolic conservation laws. Since mulitderivative schemes require more function evaluations, but smaller memory footprints to achieve greater accuracy, they fit perfectly into the computing paradigm of GPUs.

In this article, we obtain A-stable schemes of arbitrary order, using a MOL$^T$ formulation of the wave equation, and implicitly including higher order derivatives. However, we construct the derivatives using a novel approach:  recursive convolution. The cornerstone of our method lies in the fact that, using recursive applications of the convolution operators introduced in \cite{Causley2013a}, the inversions of higher derivatives can be performed analytically, so that the resulting scheme is made explicit, even at the semi-discrete level. In constructing the analytical convolution operators, we incorporate the boundary conditions directly, and as a result Dirichlet and periodic boundary conditions can be implemented to higher order with no additional complexity.

Furthermore, the convolution operator can be applied in $O(N)$ operations, and so the schemes we find in $d = 1,2$ and 3 dimensions will achieve accuracy $2P$ in $O(P^d N)$ operations per time step. Finally, the convolution algorithm utilizes exponential recurrence relations, which effectively localize contributions to the spatial integrals, making it suitable for domain decomposition. Thus, our algorithm will scale to multiple processors much more efficiently than traditional ADI solvers, which utilize global, algebraic solvers. A parallel implementation of our algorithm is the subject of future work.

The rest of this paper is laid out as follows. In Section \ref{sec:Background}, we briefly describe the main features of the MOL$^T$ algorithm. In Section \ref{sec:Scheme1d}, we derive a family of schemes of order $2P$, and prove their order of accuracy, and unconditional stability. In Section \ref{sec:SchemeADI}, we generalize the first order results to higher spatial dimensions, producing ADI methods of order $2P$ which will be A-stable. We conclude with a brief discussion in Section \ref{sec:Conclusion}.

\section{Background and notation}
\label{sec:Background}
We begin with a review of the relevant details of our method. More details can be found in \cite{Causley2013a}, \cite{Causley2013}. Consider the one-dimensional wave equation
\begin{equation}
	\label{eqn:wave}
	\frac{1}{c^2 }u_{tt} =u_{xx}, \quad x \in \mathbb{R}
\end{equation}
with prescribed initial conditions.

Using the method of lines transpose (MOL$^T$), we introduce the semi-discrete solution $u^{n} =u^{n}(x)$, which approximates $u(x,t)$ at $t = t_n = n\Delta t$. We then replace the second time derivative with the standard second order finite difference stencil, so that
\begin{equation}
	\label{eqn:Semi}
	\frac{u^{n+1}-2u^n+u^{n-1}}{(c\Delta t)^2} = \frac{\partial^2}{\partial x^2} \left(u^n+\frac{u^{n+1}-2u^n+u^{n-1}}{\beta^2}\right),
\end{equation}
where we also introduce the averaging parameter $\beta>0$, to ensure that the second spatial derivative appears implicitly, and also that the scheme remains symmetric about $t^n$.

After some rearranging of terms, we arrive at the modified Helmholtz equation, which can be written in the form
\begin{equation}
	\label{eqn:L_eq}
	\mathcal{L}\left[  u^{n+1}-2u^n+u^{n-1} \right] = \beta^2\left( u^n - \mathcal{L}\left[ u^n\right]  \right)
\end{equation}
where
\begin{equation}
	\label{eqn:L_def}
	\mathcal{L}: = 1- \frac{1}{\alpha^2}\frac{\partial^2}{\partial x^2} , \quad \alpha = \frac{\beta}{c\Delta t}, \quad 0<\beta \leq 2.
\end{equation}
The differential equation can be solved by convolution with the free space Green's function, which In 1d means that
\begin{equation}
	\label{eqn:L_Inverse}
	\mathcal{L}^{-1}[u(x)]:= \frac{\alpha}{2}\int_{-\infty}^\infty e^{-\alpha|x-y|}u(y)dy.
\end{equation}
The definition can additionally be modified to include boundary corrections on a finite domain (see \cite{Causley2013a}). We also introduce a new operator related to \eqref{eqn:L_Inverse}
\begin{equation}
	\label{eqn:D_def_1d}
	u^{n+1}-2u^n+u^{n-1} = -\beta^2\mathcal{D}[u^n], \quad \mathcal{D}[u](x):= u(x) - \mathcal{L}^{-1}[u](x),
\end{equation}
which will be used extensively in the ensuing discussion. The semi-discrete solution \eqref{eqn:D_def_1d} is therefore defined in terms of a convolution integral, and as mentioned, traditional methods of discretization in space will bear a cost of $O(N^2)$ operations per time step to evaluate $u^{n+1}(x)$ at $N$ spatial points. However, we have developed a fast convolution algorithm for \eqref{eqn:D_def_1d}, so that the numerical solution is obtained in $O(N)$ operations per time step. This is accomplished by first performing a "characteristic" decomposition $\mathcal{L}^{-1}[u](x) = I^L(x) + I^R(x)$, where
\begin{align*}
	I^L(x) = \frac{\alpha}{2}\int_{-\infty}^x u(y)e^{-\alpha(x-y)}dy, \quad I^R(x) = \frac{\alpha}{2}\int_x^{\infty} u(y)e^{-\alpha(y-x)}dy,
\end{align*}
so that both integrands decay exponentially away from $x$. Additionally, they satisfy exponential recurrence relations, which means that
\[
	I^L(x) = e^{-\alpha \delta}I^L(x-\delta) + \frac{\alpha}{2}\int_0^\delta e^{-\alpha y} u(x-y)dy, \quad 	I^R(x) = e^{-\alpha \delta}I^R(x+\delta) + \frac{\alpha}{2}\int_0^\delta e^{-\alpha y} u(x+y)dy.
\]
These expressions are exact, and upon discretization, the integrals are approximated with $O(1)$ operations at each of the $N$ points, hence resulting in an $O(N)$ scheme. We have also proven that the resulting fully discrete solution is second order accurate in time and space, and unconditionally stable (i.e., A-stable) for $0<\beta\leq 2$ \cite{Causley2013a}.

\begin{remark}
While it is not obvious from the update scheme \eqref{eqn:D_def_1d}, the solution $u^{n+1}(x)$ is the solution of an implicit scheme. This is more apparent from viewing equation \eqref{eqn:L_eq}, where $\mathcal{L}$ is acting on the unknown solution $u^{n+1}$. The fact that $u^{n+1}$ is given explicitly by \eqref{eqn:D_def_1d} is a feature of the MOL$^T$ formulation, which provides a means to analytically invert the semi-discrete Helmholtz operator $\mathcal{L}$. This is in contrast to the MOL formulation, which inverts an approximate (algebraic) spatial operator.
\end{remark}

\section{An A-stable family of schemes of order 2p}
\label{sec:Scheme1d}
As mentioned in our introductory remarks, we can achieve a higher order scheme by including more spatial derivatives in the numerical scheme. We shall continue to perform our analysis at the semi-discrete level, and make comments about the spatial discretization where appropriate. To motivate our discussion, let's first apply the Lax-Wendroff procedure to the semi-discrete wave equation, exchanging time derivatives for spatial derivatives in the Taylor expansion
\begin{align}
	\label{eqn:exact}
	u^{n+1}-2u^{n}+u^{n-1}	&= 2\sum_{m=1}^{\infty} \frac{\Delta t^{2m}}{(2m)!}\left(\partial_{tt}\right)^{m}u^n \\
						&= 2\sum_{m=1}^{\infty} \frac{\beta^{2m}}{(2m)!}\left(\frac{c\Delta t}{\beta}\right)^{2m}\left(\partial_{xx}\right)^{m}u^n \nonumber \\
						&= 2\sum_{m=1}^{\infty} \frac{\beta^{2m}}{(2m)!} \left(\frac{\partial_{xx}}{\alpha^2}\right)^{m}u^n \nonumber.
\end{align}
In the second step of this expansion, we have used the fact that
\[
	\left(\partial_{tt}\right)^m u = \left(c^2\partial_{xx}\right)^m u, \quad m \geq 1.
\]
Our next goal is to approximate the differential operators $(\partial_{xx})^m$ using the compositions of the convolution operator $\mathcal{D}$ from \eqref{eqn:D_def_1d}. We begin this process by observing
\[
	\mathcal{F}\left[\left(\frac{\partial_{xx}}{\alpha^2}\right)^{m}\right] =  (-1)^m\left(\frac{k}{\alpha}\right)^{2m}
\]
and 
\[
	\hat{D} := \mathcal{F}\left[ \mathcal{D} \right] = 1-\mathcal{F}\left[\mathcal{L}^{-1}\right] = 1 - \frac{1}{1+\left(\frac{k}{\alpha}\right)^2} = \frac{\left(\frac{k}{\alpha}\right)^2}{1+\left(\frac{k}{\alpha}\right)^2}.
\]
From this final expression for $\hat{D}$ we solve for the quantity $(k/\alpha)^2$, finding
\begin{align*}
	&\left(\frac{k}{\alpha}\right)^2 = \frac{\hat{D}}{1-\hat{D}} = \sum_{p=1}^\infty \hat{D}^{p} \quad \text{and}\\
	&\left(\frac{k}{\alpha}\right)^{2m} = \left(\frac{\hat{D}}{1-\hat{D}}\right)^{m} =  \sum_{p=m}^\infty \binom{p-1}{m-1}\hat{D}^{p},
\end{align*}
which now gives an exact expression for all even derivatives, defined solely in terms of $\mathcal{D}$. Inserting these into the Taylor expansion \eqref{eqn:exact}
\[
	u^{n+1}-2u^{n}+u^{n-1}	= \sum_{m=1}^{\infty} (-1)^m\frac{2\beta^{2m}}{(2m)!} \sum_{p=m}^\infty \binom{p-1}{m-1}\mathcal{D}^{p}[u^n].
\]
While this expression is interesting from a theoretical standpoint, it holds little appeal in practice. However, if we reverse the order of summation, then the inner sum can be collapsed, and we have
\begin{equation}
	\label{eqn:update_P}
	u^{n+1}-2u^{n}+u^{n-1}	= \sum_{p=1}^{\infty} A_p(\beta)\mathcal{D}^{p}[u^n],
\end{equation}
where the coefficients are polynomials in successively higher orders of $\beta^2$, and are given by
\begin{equation}
	A_p(\beta) = 2\sum_{m=1}^p (-1)^m\frac{\beta^{2m} }{(2m)!}\binom{p-1}{m-1}.
\end{equation}
Since
\[
	\mathcal{D}^p[u^n] \approx \left(-\frac{\partial_{xx}}{\alpha^2}\right)^p u^n + O\left(\frac{1}{\alpha^{2p+2}}\right)
\]
and $\alpha^{-1} = O(\Delta t)$, the series converges to a solution of the wave equation. Furthermore, once truncated, the $P$-term approximation will be accurate to order $2P$. For clarity, the first three schemes are
\begin{align}
	\label{eqn:update_1}
	u^{n+1}-2u^{n}+u^{n-1}	&= -\beta^2 \mathcal{D}[u^n] \\
	\label{eqn:update_2}
	u^{n+1}-2u^{n}+u^{n-1}	&= -\beta^2 \mathcal{D}[u^n]-\left(\beta^2-\frac{\beta^4}{12}\right)\mathcal{D}^2[u^n] \\
	\label{eqn:update_3}
	u^{n+1}-2u^{n}+u^{n-1}	&= -\beta^2 \mathcal{D}[u^n]-\left(\beta^2-\frac{\beta^4}{12}\right)\mathcal{D}^2[u^n]  - \left(\beta^2-\frac{\beta^4}{6}+\frac{\beta^6}{360}\right)\mathcal{D}^3[u^n].
\end{align}
The application of the operator $\mathcal{D}^p$ will be computed recursively, and so a scheme of order $P$ will have a cost of $O(PN)$ per time step for $N$ spatial discretization points. Notice that equation \eqref{eqn:update_1} is in fact identical to the original second order scheme \eqref{eqn:D_def_1d}, while the schemes \eqref{eqn:update_2} and \eqref{eqn:update_3} are fourth and sixth order, respectively. The local truncation errors will be $O((c\Delta t/\beta)^{2P+2})$, which means that the error constant will decrease with increasing $\beta$.

\subsection{Stability}
We now prove that the schemes of order $2P$ given by truncating \eqref{eqn:update_P} are unconditionally stable, for some range of the parameter $\beta$. Our main result is that, while the truncation error will decrease with increasing $\beta$, there is a maximal value $\beta_{max}$, which depends on $P$, for which the scheme remains stable. In this respect, the value $\beta_{max}$ is the optimal choice for the $P$th scheme.

We shall prove stability in the free-space case, using Von-Neumann analysis. The case of a bounded domain is similar, but consideration of specific boundary conditions must be undertaken, and thus handled individually (the Dirichlet case was shown in \cite{Causley2013} for the second order method). Upon taking the Fourier transform in space, and introducing the amplification factor
\[
	\hat{u}^n = \rho^n \hat{u}^0
\]
the scheme will be A-stable provided that $|\rho|\leq 1$. Substitution into \eqref{eqn:update_P}, and after cancellation of the common terms, we form the characteristic polynomial satisfied by the amplification factor
\begin{equation}
	\label{eqn:P_Stab}
	\rho^2-2\rho+1 = \rho S(\beta,\hat{D}),
\end{equation}
where
\begin{equation}
	\label{eqn:Stab_sum}
	S(\beta,\hat{D}) = -\left(\sum_{p=1}^P A_p(\beta) \hat{D}^p\right), \quad \beta>0, \quad 0\leq \hat{D} \leq 1.
\end{equation}

Upon applying the Cohn-Schur criterion \cite{Strikwerda1989}, we find that the scheme will be A-stable, provided that
\[
	0\leq S(\beta,\hat{D}) \leq 4.
\]
We proceed to analyze this inequality by first proving that $S$ is strictly increasing as a function of $\hat{D}$ for some interval $0<\beta \leq \beta^*$, and then by finding a maximal value $\beta_{max}$ for which stability of the scheme is ensured for any $\Delta t$. We will make use of the following 
\begin{lemma}
For each $P\geq 1$, there exists $\beta^* >0$ such that for $0<\beta \leq \beta^*$, $A_p<0$ for each $1 \leq p \leq P$.
\end{lemma}

{\em Proof}. Our proof is by induction. The case $p=1$, is trivial, since $A_1 = -\beta^2$, and so $A_1<0$ for any choice $\beta^*>0$. For $p>1$, we first choose $\beta^*$ for which $A_p$ is strictly negative, and then show that the same is automatically true for $A_k$, for all $k<p$. To choose $\beta^*$, first note that
\[
	A_p = -\frac{2}{0!2!}\beta^2\left(1-\frac{p-1}{1\cdot 3 \cdot 4}\beta^2\right)-\frac{2(p-1)(p-2)}{2! 6!}\beta^6\left(1-\frac{p-3}{3\cdot 7\cdot 8}\beta^2\right)-\ldots
\]
Thus, whenever $\beta\leq \beta^*= \sqrt{12/(p-1)}$, the first term inside the parentheses is strictly positive, and in fact all remaining terms will also be positive. Now, notice that $\beta^*$ will decrease monotonically with increasing $p$, and so if we choose $\beta^*$ to ensure that $A_p<0$, then it immediately follows that $A_k<0$ for all $k\leq p$.
\endproof

We are now prepared to state our main
\begin{theorem}
The semi-discrete scheme, given by truncating the sum in \eqref{eqn:update_P} after $P$ terms, will be unconditionally stable, provided that $0< \beta \leq \beta_{max}$, where
\begin{equation}
	\label{eqn:Stab_cond}
	-\sum_{p=1}^P A_p(\beta_{max}) = 4
\end{equation}
\end{theorem}

{\em Proof}.
As a result of the lemma, we are guaranteed an interval $(0,\beta^*)$ which, for fixed $P$, all $A_p<0$ for $1\leq p \leq P$. Therefore, the sum \eqref{eqn:Stab_sum} is strictly positive, and increasing in both $\beta$ and $\hat{D}$. Thus, we only need to study the extremal value $\hat{D} \to 1$, and $\beta_{max}$, which we find by solving the equality \eqref{eqn:Stab_cond}.
\endproof

\begin{figure}[htb]
\begin{centering}
	\includegraphics[width=.45\textwidth]{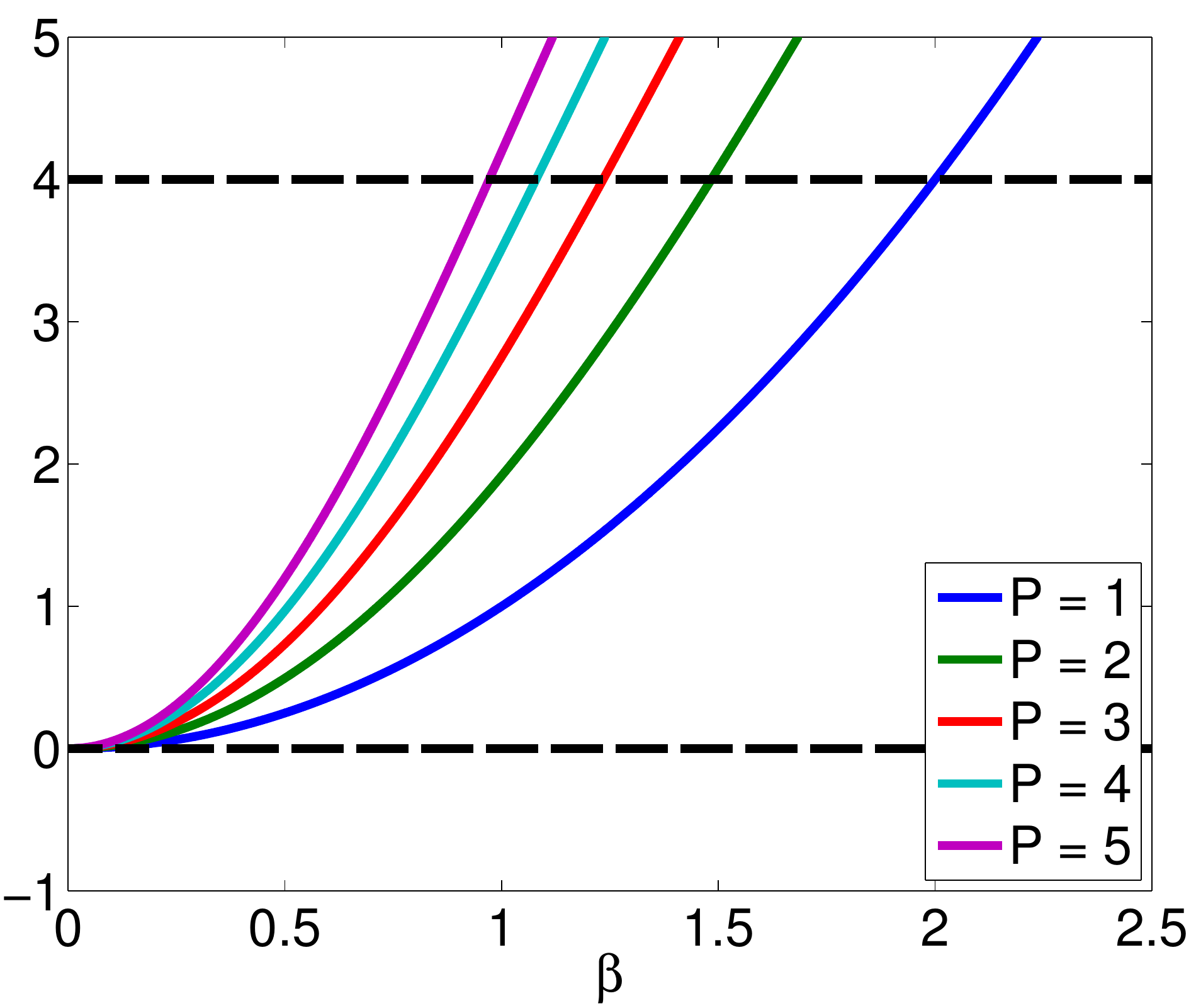}
	\caption{A plot of the the sum \eqref{eqn:Stab_sum}, corresponding to the scheme \eqref{eqn:update_P} truncated at $P$ terms. The scheme is A-stable whenever $S(\beta,\hat{D})\leq 4$.}
	\label{fig:Stab_plot}
\end{centering}
\end{figure}
Hence the stability condition amounts to finding the first positive root of a polynomial of order $P$ in $\beta^2$. A plot of the stability polynomial \eqref{eqn:Stab_sum} with $\hat{D}=1$ is shown for the first few values of $P$ in Figure \ref{fig:Stab_plot}. The value $\beta_{max}$ is taken as the intersection of the curve with the dotted line $S = 4$. The values for these first few schemes are also shown in Table \ref{tab:Stab_table}.
\begin{remark}
While $\beta_{max}$ decreases as $P$ increases, the double root at $\beta=0$ will guarantee the existence of some region $(0,\beta^*)$ for which stability can be achieved. Additionally, methods which include more than $P$ terms can be derived, which achieve only order $2P$, but allow for produce a smaller error constant by increasing $\beta_{max}$. A rigorous investigation of their construction has not yet been undertaken.
\end{remark}

\begin{table}
\begin{center}
\begin{tabular}{ | l  ||c|c|c|c|c |}
\hline
	P			& 1	& 2		& 3		& 4		& 5	\\
\hline
	Order		& 2	& 4		& 6		& 8		& 10	\\
\hline
	$\beta_{max}$	& 2	& 1.4840	& 1.2345	& 1.0795	& 0.9715 \\
\hline
\end{tabular}
\caption{The maximum values $\beta_{max}$ for which the $P$-th scheme remains A-stable, is found by solving $S(\beta_{max},1) = 4$.}
\label{tab:Stab_table}
\end{center}
\end{table}

\subsection{High order initial and boundary conditions}
Since the scheme \eqref{eqn:update_P} is a 3-step method, it will require two initial starting values, which must be computed to $O(\Delta t^{2P})$ in order for the numerical solution to achieve the expected order. While the initial condition $u^0 = f(x)$ is imposed exactly, the value $u^{1} = u(x,\Delta t)$ must be approximated. In analogy to the derivation above, we shall proceed with a Taylor expansion, and using the Lax-Wendroff procedure, convert all even time derivatives into spatial derivatives. However, the odd time derivatives will instead make use of the initial velocity, $u_{t}(x,0) = g(x)$. Thus
\begin{align}
	u^1	&= \sum_{m=0}^\infty \frac{\Delta t^m}{m!}\partial_t^mu^0 \nonumber \\
		&= \sum_{m=0}^\infty \left(\frac{\Delta t^{2m}}{(2m)!}\partial_ t^{2m}u^0+ \frac{\Delta t^{2m+1}}{(2m+1)!}\partial_t^{2m+1}u^0\right) \nonumber \\
		&= \sum_{m=0}^\infty \left(\frac{(c\Delta t)^{2m}}{(2m)!}\partial_ x^{2m}u^0+ \frac{(c\Delta t)^{2m+1}}{(2m+1)!}\partial_x^{2m}\frac{1}{c}\partial_t u^0\right) \nonumber \\
		&= \sum_{m=0}^\infty \left(\frac{(c\Delta t)^{2m}}{(2m)!}\partial_ x^{2m}f(x)+ \frac{(c\Delta t)^{2m+1}}{(2m+1)!}\partial_x^{2m}\frac{1}{c}g(x)\right).
\end{align}
This expansion can now be truncated at $O(\Delta t^{2P})$, and since all spatial derivatives are even, they can be approximated using convolution \eqref{eqn:D_def_1d}. Alternatively they can be computed analytically, since $f$ and $g$ are known.

We also consider the application of boundary conditions for recursive applications of the convolution operator \eqref{eqn:D_def_1d}. We demonstrate the approach for Dirichlet boundary conditions; suppose $u(a,t) = u_L(t)$, and $u(b,t) = u_R(t)$ are prescribed. Since $\mathcal{D}^m[u] \approx (\partial_{xx}/\alpha^2)^m u$, we require even spatial derivatives of $u$ at $x = a$ and $b$. But, since $u$ satisfies the wave equation for $a<x<b$, we can use the inverse Lax-Wendroff procedure, and upon taking the appropriate limit, find
\begin{equation}
	\lim_{x\to a}\left(\partial_{xx}\right)^m  u = \lim_{x \to a} \left(\frac{\partial_{tt}}{c^2}\right)^m u = \left(\frac{\partial_{tt}}{c^2}\right)^m  u_L(t),
\end{equation}
with the corresponding result for holding for $x=b$. The result for Neumann boundary conditions will also follow from a similar procedure, but upon considering odd derivatives. Additionally, periodic boundary conditions can be implemented in a straightforward manner, as was shown in the basic algorithm in \cite{Causley2013a}.

\subsection{Numerical results}
\begin{figure}[htb]
\begin{centering}
	\includegraphics[width=.45\textwidth]{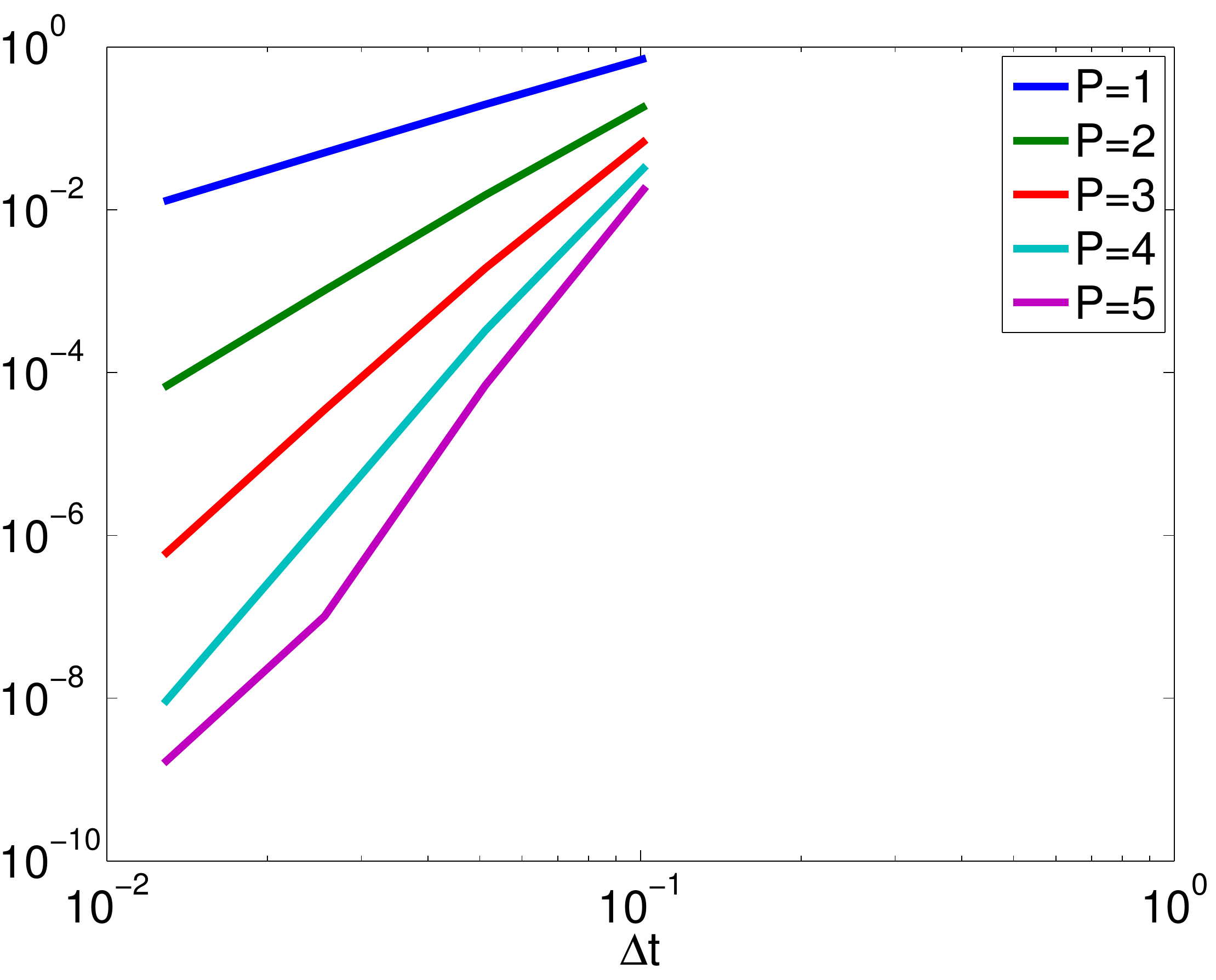}
	\caption{Convergence in the $L^2$-norm of a 1d standing wave, with Dirichlet boundary conditions. The spatial resolution is held fixed at $\Delta x = 0.0001$.}
	\label{fig:Conv_plot}
\end{centering}
\end{figure}

Before moving onward to multi-dimensional schemes, we first illustrate the accuracy of our method for a 1d example. We perform time marching for a standing wave $u(x,0) = \sin(2\pi x)$, for $x \in [0,1]$, up to time $T = 1$. Since the fast convolution algorithm from \cite{Causley2013a} is second order accurate in space, we fix $\Delta x = 0.0001$ to ensure that the dominant error in the solution is temporal. The $L^2$-norm of the error is plotted in Figure \ref{fig:Conv_plot} for several values of $\Delta t$, with varying order $P$. For each $P$, we used $\beta = \beta_{max}$ according to Table \eqref{tab:Stab_table}. In the 10th order scheme, the spatial error can be seen to dominate the error for the smallest value $\Delta t$. Thus, refining further in time would produce no further improvement.

\section{Higher dimensions and higher accuracy via ADI splitting}
\label{sec:SchemeADI}
We next address the solution of the wave equation \eqref{eqn:wave} in higher dimensions using alternate direction implicit (ADI) splitting. To achieve schemes of higher order, we again  begin with the expansion \eqref{eqn:exact}, but now the Lax procedure introduces the Laplacian
\[
	u^{n+1}-2u^n+u^{n-1}= 2\sum_{m=1}^\infty \frac{\Delta t^{2m}}{(2m)!}\left(\frac{\partial^{2m}}{\partial t^{2m}}\right)u^n= 2\sum_{m=1}^\infty \frac{\beta^{2m}}{(2m)!}\left(\frac{\nabla^2}{\alpha^2}\right)^m u^n.
\]
In order to approximate higher order powers of the Laplacian operator using ADI splitting, we first define univariate modified Helmholtz operators, and their corresponding $\mathcal{D}$ operators as
\begin{equation}
	\label{eqn:LD_gam}
	\mathcal{L}_\gamma:= 1-\frac{\partial_\gamma^2}{\alpha^2}, \quad \mathcal{D}_\gamma := 1 - \mathcal{L}_\gamma^{-1} , \quad \gamma = x,y,z.
\end{equation}
Notice that these operators satisfy the following identity
\begin{equation}
	\label{eqn:LD_identity}
	\mathcal{L}_\gamma\mathcal{D}_\gamma[u] = \mathcal{L}[u] - u = -\frac{\partial_{\gamma\gamma}}{\alpha^2}u, \quad \gamma = x,y,z.
\end{equation}
Thus the Laplacian is given by
\[
	-\frac{\nabla^2}{\alpha^2} = \mathcal{L}_x\mathcal{D}_x+\mathcal{L}_y\mathcal{D}_y+\mathcal{L}_z\mathcal{D}_z = \mathcal{L}_x\mathcal{L}_y\mathcal{L}_z [\mathcal{C}_{xyz}],
\]
where the new operator is
\begin{equation}
	\label{eqn:C_def_ADI_3}
	\mathcal{C}_{xyz} := \mathcal{L}_y^{-1}\mathcal{L}_z^{-1}\mathcal{D}_x+\mathcal{L}_z^{-1}\mathcal{L}_x^{-1}\mathcal{D}_y+\mathcal{L}_x^{-1}\mathcal{L}_y^{-1}\mathcal{D}_z.
\end{equation}
The corresponding 2d operator is
\begin{equation}
	\label{eqn:C_def_ADI}
	\mathcal{C}_{xy} := \mathcal{L}_y^{-1}\mathcal{D}_x+\mathcal{L}_x^{-1}\mathcal{D}_y.
\end{equation}
This result, while interesting, is not quite satisfactory. The first issue is that the Laplacian is now given in terms of $\mathcal{L}_\gamma$, and thus requires approximations of spatial derivatives, which we are trying to avoid. Secondly, the directional sweeps of the ADI convolutions are represented by a composition of the operators $\mathcal{L}_\gamma^{-1}$, not a sum of them. In two dimensions, the ADI scheme is defined by
\begin{equation}
	\label{eqn:D_def_ADI}
	\mathcal{D}_{xy} := 1-\mathcal{L}_x^{-1}\mathcal{L}_y^{-1},
\end{equation}
while in three dimensions it is
\begin{equation}
	\label{eqn:D_def_ADI_3}
	\mathcal{D}_{xyz} := 1-\mathcal{L}_x^{-1}\mathcal{L}_y^{-1}\mathcal{L}_z^{-1}.
\end{equation}
Motivated by this form we appeal to one final identity, obtained by rearranging and inverting \eqref{eqn:D_def_ADI_3},
\begin{align*}
	\mathcal{L}_x\mathcal{L}_y\mathcal{L}_z =\left(1-\mathcal{D}_{xyz}\right)^{-1},
\end{align*}
which means that
\[
	-\frac{\nabla^2}{\alpha^2} =\left(1-\mathcal{D}_{xyz}\right)^{-1} \mathcal{C}_{xyz}.
\]
Thus, all even order of the Laplacian can be constructed by expanding the symbol $(1-\mathcal{D})^{-m}$ as a power series (which is now in terms of the ADI operator!), and the result is
\begin{equation}
	\label{eqn:Laplacian_2D_ADI_3}
	\left(\frac{\nabla^2}{\alpha^2}\right)^m
	 =  (-1)^m\mathcal{C}_{xyz}^m \sum_{p=m}^\infty \binom{p-1}{m-1}\mathcal{D}_{xyz}^{p-m}.
\end{equation}
The corresponding 2d result is
\begin{equation}
	\label{eqn:Laplacian_2D_ADI}
	\left(\frac{\nabla^2}{\alpha^2}\right)^m
	 =  (-1)^m\mathcal{C}_{xy}^m \sum_{p=m}^\infty \binom{p-1}{m-1}\mathcal{D}_{xy}^{p-m},
\end{equation}
which can be seen to have the identical form in this notation, except for the spatial subscripts. Upon omitting the subscripts, the semi-discrete scheme for 2d and 3d is
\begin{align}
	u^{n+1}-2u^n+u^{n-1}	&= \sum_{m=1}^\infty \frac{2\beta^{2m}}{(2m)!}\left(\frac{\nabla^2}{\alpha^2}\right)^m u^n
						\nonumber \\
						&= \sum_{m=1}^\infty (-1)^m\frac{2\beta^{2m}}{(2m)!} \mathcal{C}^m \sum_{p=m}^\infty \binom{p-1}{m-1}\mathcal{D}^{p-m} [u^n]
						\nonumber \\
						\label{eqn:Full_2D}
						&= \sum_{p=1}^\infty  \sum_{m=1}^p (-1)^m\frac{2\beta^{2m}}{(2m)!}\binom{p-1}{m-1} \mathcal{C}^m\mathcal{D}^{p-m} [u^n].
\end{align}
It is interesting to compare to the 1d scheme \eqref{eqn:update_P}, which can in fact be recovered by setting $\mathcal{C} = \mathcal{D}$.

Upon truncation at $p = P$, we obtain a scheme of order $2P$, the first few of which are
\begin{align*}
	u^{n+1}-2u^n+u^{n-1}&= -\beta^2\mathcal{C}[u^n] \\
	u^{n+1}-2u^n+u^{n-1}&= -\beta^2\mathcal{C}[u^n] - \left(\beta^2\mathcal{D} - \frac{\beta^4}{12}\mathcal{C}\right)\mathcal{C}[u^n] \\
	u^{n+1}-2u^n+u^{n-1}&= -\beta^2\mathcal{C}[u^n] - \left(\beta^2\mathcal{D} - \frac{\beta^4}{12}\mathcal{C}\right)\mathcal{C}[u^n]-\left(\beta^2\mathcal{D}^2 - \frac{\beta^4}{6}\mathcal{C}\mathcal{D}+\frac{\beta^6}{360}\mathcal{C}^2\right)\mathcal{C}[u^n].
\end{align*}
As in the 1d case, these schemes will be unconditionally stable for all $\Delta t$, and the same range for $\beta$ as shown in Table \ref{tab:Stab_table}.

\subsection{Inclusion of source terms}
Until now we have considered raising the order of approximations only for the homogeneous wave equation. We now consider the inclusion of general source terms $S(u,x,t)$, and the higher order schemes generated by including higher time derivatives. To do so, we consider
\begin{equation}
	\frac{1}{c^2}u_{tt} = \nabla^2 u + S,
\end{equation}
and upon taking even order time derivatives, find
\begin{align*}
	\left(\frac{\partial_{tt}}{c^2}\right)^m u &= \left(\frac{\partial_{tt}}{c^2}\right)^{m-1} \left( \nabla^2 u + S\right) \\
	&= \left(\frac{\partial_{tt}}{c^2}\right)^{m-2}\nabla^2 \left( \nabla^2 u + S\right)+\left(\frac{\partial_{tt}}{c^2}\right)^{m-1}S \\
	&= \vdots \\
	&= \left(\nabla^2\right)^m u + \left[ \left(\frac{\partial_{tt}}{c^2}\right)^{m-1} +\nabla^2\left(\frac{\partial_{tt}}{c^2}\right)^{m-2}+\ldots \left(\nabla^2\right)^{m-1}\right] S.
\end{align*}
Thus, in order for the source to be included to high accuracy, we shall require the computation of terms of the form
\[
	\left[\sum_{k=0}^{m} \left(\nabla^2\right)^k\left(\frac{\partial_{tt}}{c^2}\right)^{m-k}\right]S, \quad m=0, 1, \ldots P-1.
\]
The method of construction of such terms depends on the specific nature of the source term, specifically whether the time derivatives are approximated using finite differences, or by replacing them with other known information (i.e. from a constitutive relation). For these reasons, we only briefly consider the fourth order accurate implementation. Upon discretization of the second time derivative, and approximation of the Laplacian using the convolution operators \eqref{eqn:C_def_ADI}, we have
\[
	u^{n+1}-2u^n+u^{n-1} = -\beta^2\mathcal{C}[u^n] - \left(\beta^2\mathcal{D} - \frac{\beta^4}{12}\mathcal{C}\right)\mathcal{C}[u^n] +\frac{\beta^2}{12\alpha^2}\left(S^{n+1}+10S^n+S^{n-1}\right)+\frac{\beta^2}{12\alpha^4}\mathcal{C}[S^n].
\]
If $S$ depends either on $u$, or another dependent variable which is coupled to $u$, then this equation must be solved iteratively now, due to the appearance of source terms at time level $t_{n+1}$.

\subsection{Efficient computation and operation count of higher order methods}
One potential pitfall to the numerical implementation of the schemes \eqref{eqn:Full_2D} is the unnecessary duplication of work in constructing additional terms, which are all defined by multiple convolutions. Here we estimate the operation count per time step for the $P$th scheme, and in doing so consider two competing forces:  the economical re-use of computed quantities, and the symmetrization to remove anisotropy.

\begin{table}[hb]
\begin{center}
\begin{tabular}{|c|c|c|c|c|}
\hline
	P&Input			&Output  & Operations & Total operations\\
\hline
	1&$u^n$			&$v_1:=\mathcal{C}[u^n]$  & $O(N)$ & $O(N)$ \\
\hline
	2&$v_1$			&$v_2:=\mathcal{D}[v_1], \quad v_1:=\mathcal{C}[v_1]$ & $2O(N)$ & $3O(N)$ \\
\hline
	3&$v_1, \quad v_2$	&$v_3:=\mathcal{D}[v_2], \quad v_2:=\mathcal{C}[v_2], \quad v_1:=\mathcal{C}[v_1]$  & $3O(N)$ & $6O(N)$ \\
\hline
\end{tabular}
\caption{Operation count based on efficient re-use of computed quantities for schemes of order $2P$.}
\label{tab:Efficiency}
\end{center}
\end{table}

The first of these considerations comes from the binomial-like structure of higher order terms in the expansions. Notice that in the $P$th scheme an additional $P$ terms appear when compared to the ($P$-1)st scheme, all of which are of the form $\mathcal{C}^m\mathcal{D}^{P-m}[u^n]$. Furthermore, they can be constructed solely in terms of the previous $P$-1 terms of the previous stage, without use of any subsequent terms. Each new term requires one application of either $\mathcal{C}$, or $\mathcal{D}$, both of which can be computed in $O(N)$ operations. Thus, we obtain a simple estimate for the complexity as $O(P^2N)$ for the operation count. The procedure is demonstrated in Table \ref{tab:Efficiency} for the first few values of $P$, which also reveals that $P$ auxiliary variables $v_p$, will also be required, in addition $u^n$ and $u^{n-1}$.

The second important consideration is a practical aspect inherent to ADI schemes, which is the introduction of numerical anisotropy. To this end, it is prudent to apply the spatial convolution operators $\mathcal{L}_\gamma$ for $\gamma = x,y,z$ and average all permutations to reduce the anisotropy. This will inherently trade accuracy for computational efficiency, and as such is not included in the efficiency estimates of Table \ref{tab:Efficiency}.

\subsection{Numerical Results}

\begin{table}[ht!]
\begin{center}
\begin{tabular}{|c|c|c|c||c|c|c||c|c|c|}
\hline
			&	\multicolumn{3}{ c|| }{$P=1$}	&\multicolumn{3}{ c|| }{$P=2$}		&\multicolumn{3}{ c| }{$P=3$} \\
\hline
	$\Delta t$	&Error		& Rate	& Time (s)	&Error		& Rate	& Time (s)	&Error		& Rate	& Time (s) \\
\hline
	$0.4$	&$7.81$E-1	& * 		& 1.5		&$7.19$E-1	& * 		& 4.1		&$8.16$E-1	& * 		& $8.2$ \\
\hline
	$0.2$	&$2.46$E-1	& 1.67	& 3.9		&$1.07$E-1	& 2.74 	& 11.2	&$7.80$E-2	& 3.39	& $19.7$ \\
\hline
	$0.1$	&$7.15$E-2	& 1.78	& 7.1		&$1.03$E-2	& 3.38	& 23.3	&$2.83$E-3	& 4.78	& $44.1$ \\
\hline
	$0.05$	&$1.89$E-2	& 1.92	& 15.1	&$7.36$E-4	& 3.81	& 48.3	&$5.74$E-5	& 5.63	& $90.0$ \\
\hline
	$0.025$	&$4.84$E-3	& 1.96	& 30.0	&$4.97$E-5	& 3.89	& 94.2	&$2.29$E-6	& 4.64	& $186.2$ \\
\hline
\end{tabular}
\caption{Refinement and computational efficiency for a 2d rectangular mode $u(x,y,0) = \sin(\pi x)\sin(\pi y)$. The mesh is held fixed at $\Delta x = \Delta y = 0.003125$.}
\label{tab:ADI}
\end{center}
\end{table}

\begin{figure}[ht!]
\begin{centering}
	\subfigure[2nd order]{\includegraphics[width=.45\textwidth]{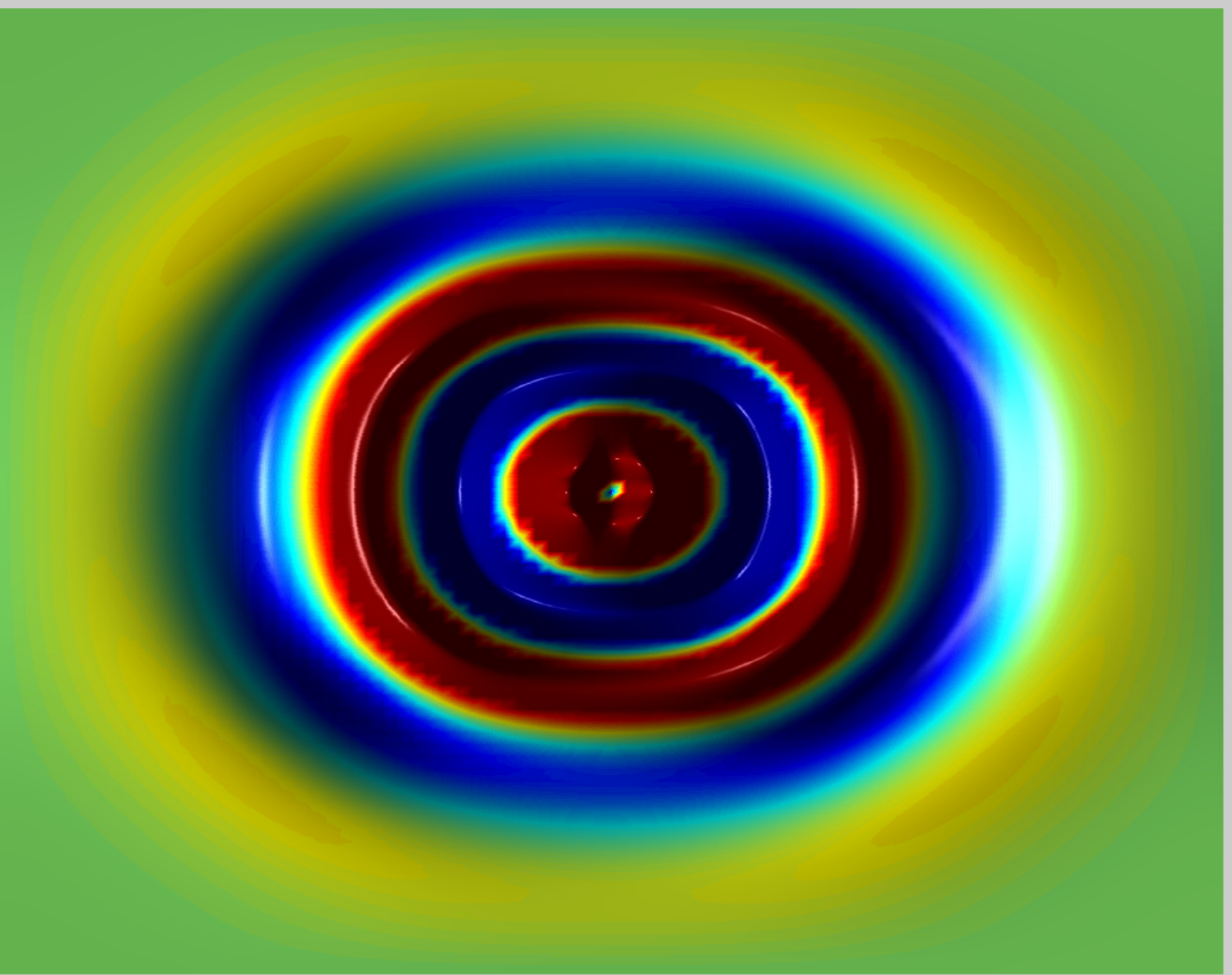}}
	\subfigure[4th order]{\includegraphics[width=.45\textwidth]{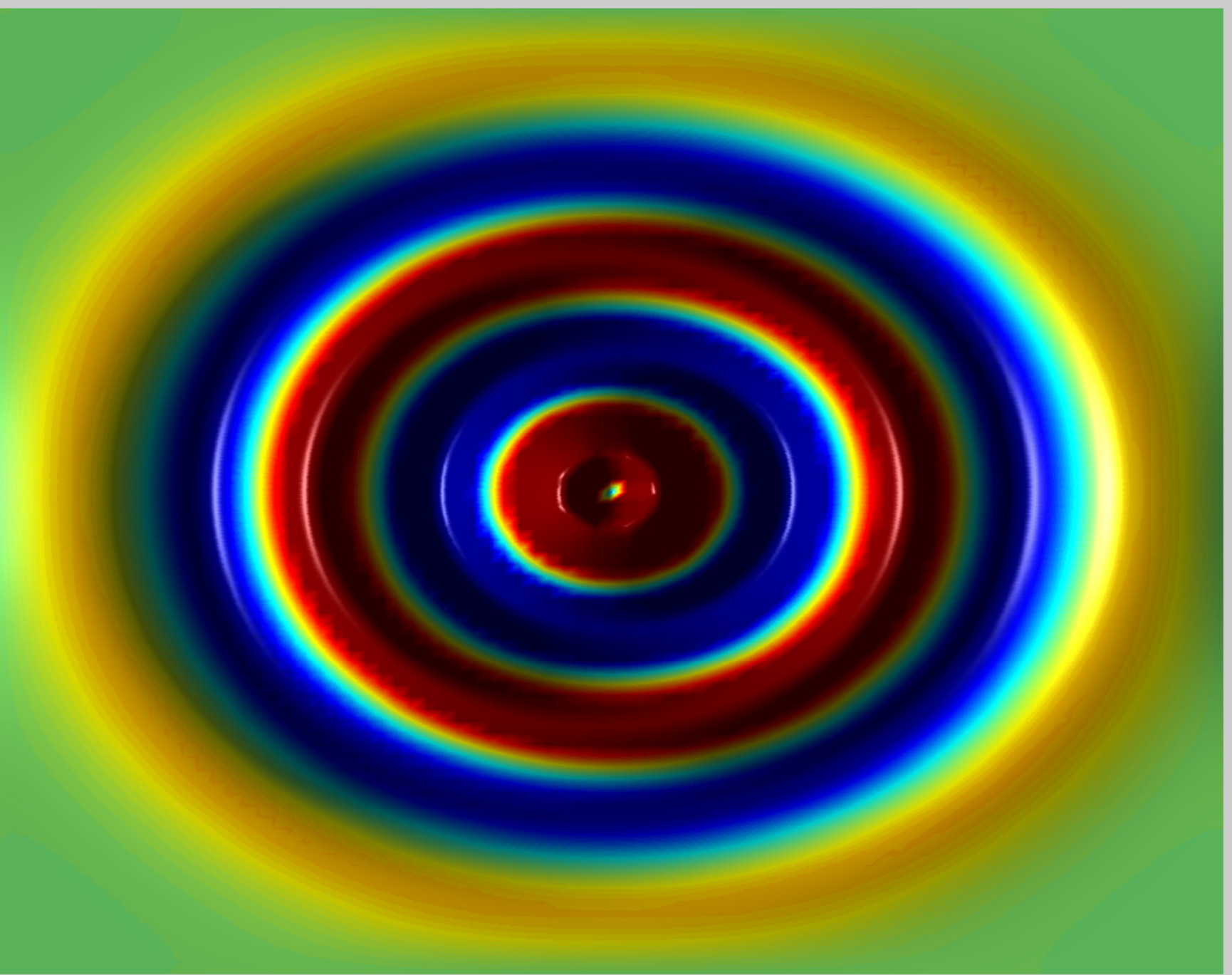}}
	\caption{Propagation due to a point source in 2d, on a $80\times 80$ mesh,with CFL number 2. The improvement due to the higher order corrections is quite apparent.}
	\label{fig:ADI}
\end{centering}
\end{figure}

We first show that our method behaves as expected, by performing a refinement study on a square domain $\Omega = [-1,1]\times[-1,1]$, with a standing mode $u(x,y,0) = \sin(\pi x)\sin(\pi y)$, up to time $T=1.2$, with a fixed spatial resolution of $641\times 641$ spatial points. The discrete $L^2$-norm of the error is constructed at each time step, and we report the maximum over all time steps in Table \ref{tab:ADI}. We additionally record the computation time required for each scheme for $P=1,2$ and 3, confirming the predicted scaling of the method from Table \ref{tab:Efficiency}.

As stated in the introduction, one major drawback of using an ADI method is the anisotropy introduced in the leading order truncation error. In Figure \ref{fig:ADI}, a sinusoidal point source located at the center of the domain is smoothly turned on, and propagated using the second order and fourth order schemes. The anisotropy is quite evident at the wavefront in the second order scheme, but is removed by implementing the fourth order scheme.

Finally, we demonstrate the advantages of a higher order method in an elliptical geometry, which is of interest in antennae design, among other applications \cite{Boriskin2008}. Currently, the most common method for simulating wave propagation in elliptical cavities is the conformal finite-difference time-domain method (CFDTD) method \cite{Dey1997}, which uses a conformal mapping to accurately represent curvilinear geometries, thus avoiding the reduction to first order due to the stair-step approximation in traditional FDTD algorithms.
\begin{figure}[hb]
\begin{centering}
	\subfigure[x-sweep]{\includegraphics[width=.43\textwidth]{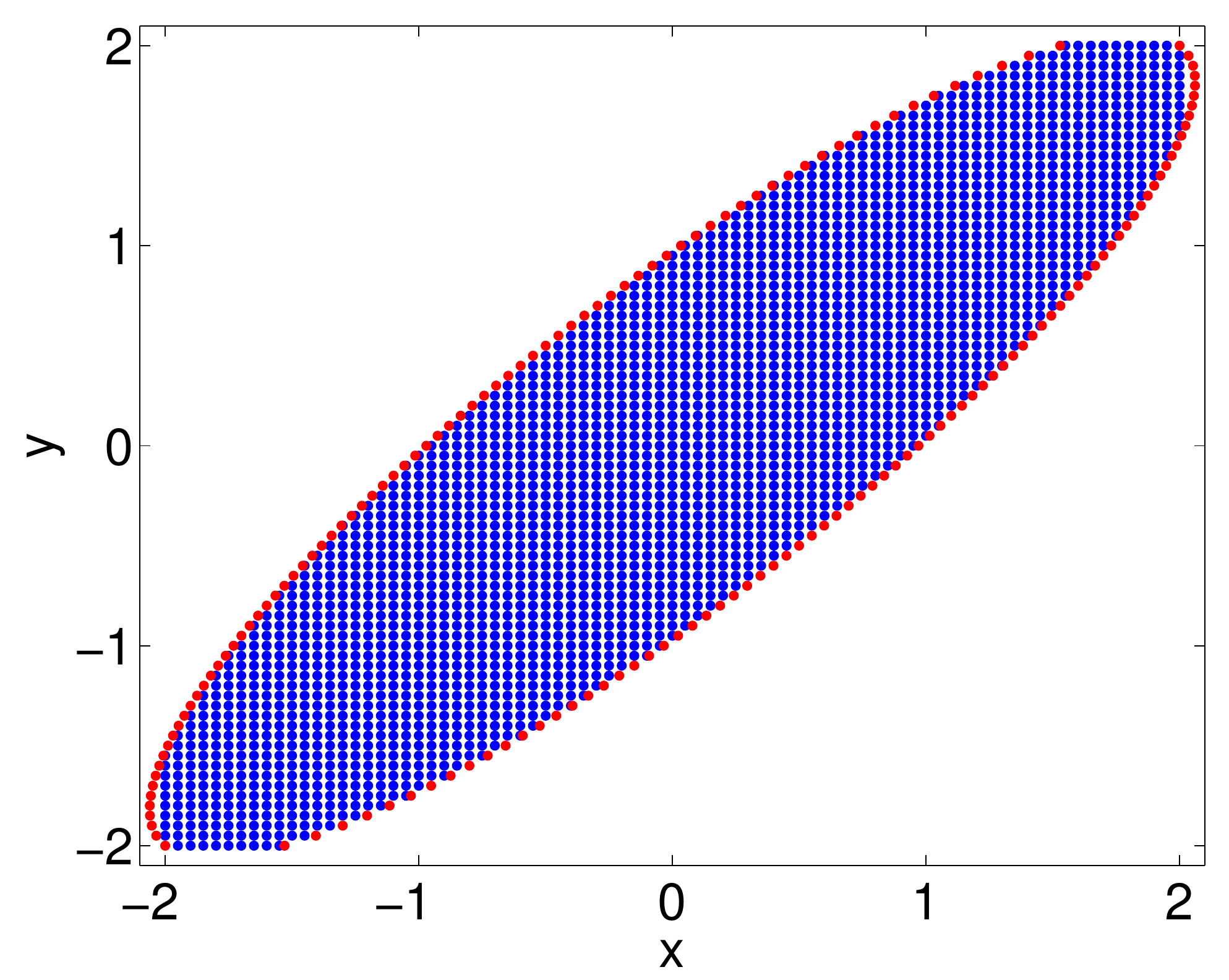}}
	\subfigure[y-sweep]{\includegraphics[width=.43\textwidth]{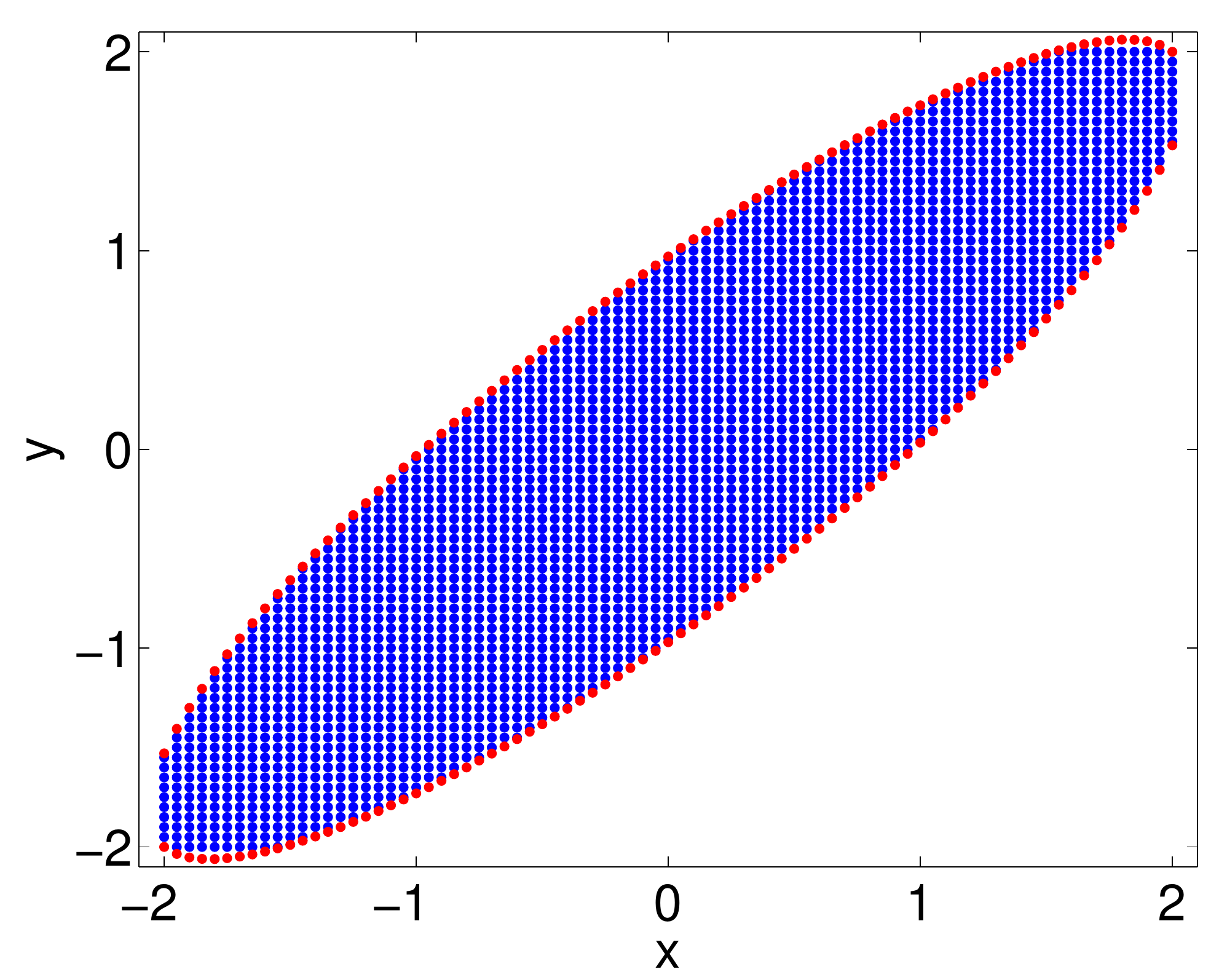}}
	\caption{Discretization of the ellipse, showing the regular Cartesian points (blue), and the additional boundary points (red) for the $x$ and $y$ sweeps.}
	\label{fig:Ellipse_points}
\end{centering}
\end{figure}

Our approach is different; rather than using conformal geometry, we embed the boundary in a regular Cartesian mesh, including a boundary point for each intersection of the ADI lines $x = x_j$, and $y = y_k$ with the boundary curve, as illustrated in Figure \ref{fig:Ellipse_points}. The $x$ and $y$ convolutions then operate on line objects, which are defined by a collection of interior points, and two boundary points, one at either end. The boundary points can be arbitrarily close to the interior points (a Taylor expansion of the coefficients is used when the spacing is small enough). See \cite{Causley2013a} for more details.

\begin{figure}[htb!]
\begin{centering}
	\subfigure{\includegraphics[width=.32\textwidth]{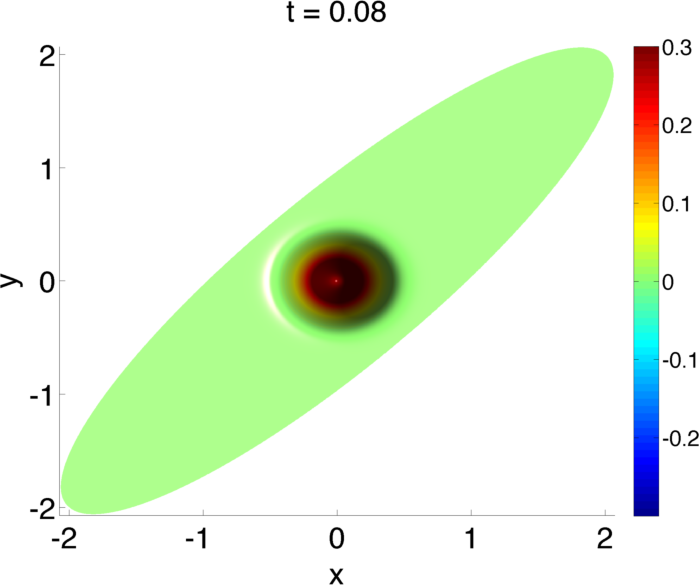}}
	\subfigure{\includegraphics[width=.32\textwidth]{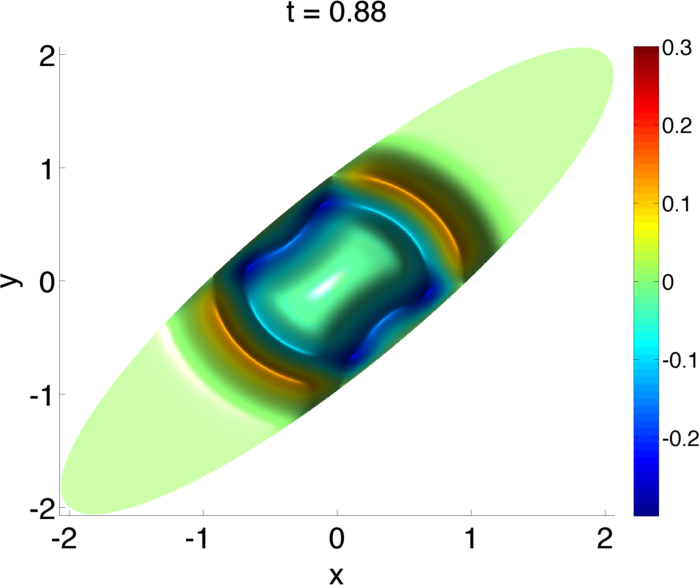}}
	\subfigure{\includegraphics[width=.32\textwidth]{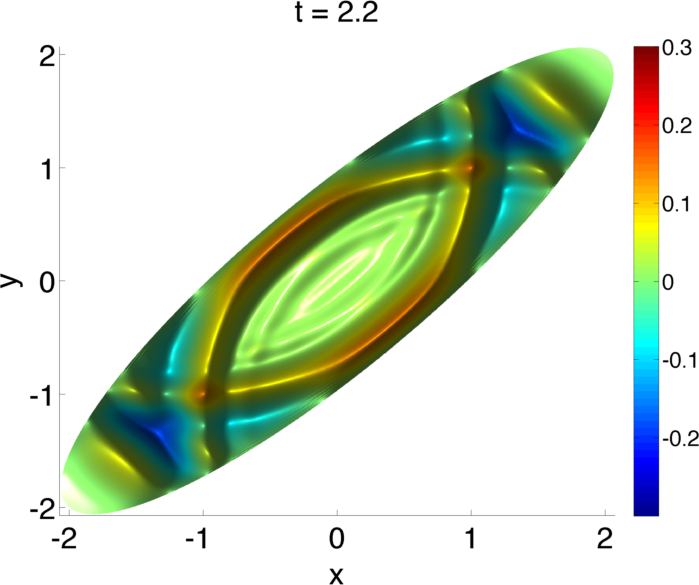}}
	\caption{Time evolution of a Gaussian field through an elliptical cavity.}
	\label{fig:Ellipse}
\end{centering}
\end{figure}

In Figure \ref{fig:Ellipse}, a Gaussian initial condition is placed inside the ellipse, whose boundary is given by
\[
	C = \left\{ (x,y): \left(\frac{x+y}{4}\right)^2+(x-y)^2 =1 \right\}.
\]
The scattering will therefore be a superposition of the natural modes of the ellipse, which are Mathieu functions. This is a great test of the algorithm, not only due to the curved boundaries, but also because the principal axes of the ellipse do not coincide with the horizontal and vertical ADI lines.

\section{Conclusions}
\label{sec:Conclusion}
In this paper we have proposed a family of schemes for the wave equation, of order $2P$, which remain A-stable by using a multi-derivative scheme in time. To maintain efficiency, we utilize a Lax approach to replace even order time derivatives with the powers of the Laplacian, which is then constructed using recursive applications of fast convolution operators previously developed for the base scheme in \cite{Causley2013a}. The resulting schemes therefore scale as $O(P^dN)$ where $d$ is the number of spatial dimensions. The expected algorithmic efficiency and convergence properties have been demonstrated in with 1d and 2d examples. This method holds great promise for developing higher order, parallelizable algorithms for solving hyperbolic PDEs, and can also be extended to parabolic PDEs.

\section*{Acknowledgements}
This work has been supported in part by AFOSR grants FA9550-11-1-0281, FA9550-12-1-0343 and FA9550-12-1-0455, NSF grant DMS-1115709, and MSU Foundation grant SPG-RG100059.


\bibliographystyle{amsplain}
\bibliography{MOLT_3.bib}


\end{document}